\documentclass{article}
\usepackage{amsmath,amsthm,amsfonts,amssymb}
\usepackage{standalone}
\usepackage{float}
\usepackage{pgfplots}
\usepackage{verbatim}
\pgfplotsset{compat=1.15}
\usepackage{tikz}
\usetikzlibrary{cd}
\usepackage{authblk}
\usepackage{mathrsfs}

\usepackage[all]{xy}
\usepackage{mathrsfs}
\usetikzlibrary{arrows}

\definecolor{uuuuuu}{rgb}{0.26666666666666666,0.26666666666666666,0.26666666666666666}

\newtheorem*{theorem*}{Theorem}
\newtheorem{theorem}{Theorem}

\theoremstyle{definition}

\newtheorem{remark}[theorem]{Remark}

\providecommand{\keywords}[1]
{
  \small	
  \textbf{\textit{Keywords---}} #1
} 

\providecommand{\AMSclass}[1]
{
  \small	
  \textbf{{2020 Mathematics subject classification---}} #1
} 

\begin{document}
\title{Ruler and compass constructions in the Lemniscate and the 17-gon}


\author[1]{M.A. G\'omez-Molleda \thanks{gomezma@uma.es}}
\author[2]{Joan-C. Lario \thanks{joan.carles.lario@upc.edu}}

\affil[1]{Universidad de M\'alaga }
\affil[2]{ UPC, Barcelona  }

\date{}

\maketitle

\begin{abstract}
We present
several ruler and compass practical geometric  constructions that can be performed in the lemniscate curve. To be precise, we provide recipes for 
halving, doubling, adding, subtracting, and transferring
lemniscate arcs with ruler and compass. 
This note complements the instructions
for the lemnatomic equilateral triangle and pentagon discussed  in \cite{GoLa}, giving the details for the construction of the lemnatomic  regular $17$-gon.
\end{abstract}

 
\AMSclass{11R20}

\section{Introduction.}\label{intro}

Right after stating the 
following
theorem on the division of the lemniscate in his {\sl Recherches sur les fonctions elliptiques}
\cite{Abel} (section V),

\begin{theorem*}[Abel, 1827] The lemniscate can be divided into $n$ equal parts with ruler and compass if $n=2^ap_1 p_2 \dots p_m,$ where $p_1, \dots, p_m$ are distinct Fermat primes,
\end{theorem*}
\noindent Abel points out its similarity with Gauss's work on the circumference:
{\sl ``Ce th\'eor\`eme est, comme on voit, pr\'ecis\'ement le m\^eme que celui de M. Gauss, relativement au cercle"}.

\bigskip
Gauss's theorem, asserting that the division of the circle in $n$ equal parts is possible with ruler and compass when $n$ is the product of a power of 2 and distinct Fermat primes, appears in paragraph 365 of the Disquisitiones Arithmeticae \cite{Gauss}, at the end of Chapter VII. Along the previous paragraphs  Gauss studies the 
 equations of the form $x^n-1=0,$ whose roots are the vertices of a regular polygon of $n$ sides on the unit circle: $\zeta_n^k=\cos{2k\pi /n} +i \sin{2k\pi/n}, \, k=0,1, \dots, n-1.$
 When~$p$ is a prime number, Gauss develops his theory of {\sl periods}, constructing a sequence 
 of algebraic numbers
 $\alpha_1, \alpha_2, \dots, \alpha_s$ such that
 $$\mathbb Q=E_0 \subset E_1\subset E_2 \subset \dots \subset E_s=\mathbb Q(\zeta_p),$$
  with 
 $E_j=E_{j-1}(\alpha_j)$ for every $j=1, \dots, s,$
 is a tower of field extensions of prime degree. 
 In Gauss's words, the division of the circle in $p$ parts has been reduced to the solution of as many equations as the number of prime factors (counting with multiplicity) of $p-1,$ and their degrees are precisely these prime factors. Gauss solves these equations by radicals using Lagrange's resolvents. When $p-1$ is a power of 2, what happens for the 
  known Fermat primes,
 $$3,5,17,257,65537,$$
 the division of the circle is reduced to equations of degree 2, and the trigonometric 
 numbers $\cos{2 k \pi/n}$ can be expressed by square roots, so that the division of the circle can be determined by geometric constructions. 

 \smallskip
 The four quadratic equations for $p=17$ are calculated in paragraph 354  of \cite{Gauss}:

$$\zeta_{17}=\frac{1}{2}\alpha_3+\frac{1}{2}\sqrt{\alpha_3^2-4} \mbox{ is a root of } x^2-\alpha_3 \, x+1, \mbox{ where }$$
$$\alpha_3=\frac{1}{2} \alpha_2+\frac{1}{2}\sqrt{2\alpha_2^3+\alpha_2^2-12\alpha_2+6} \, \mbox{ is a root of } x^2-\alpha_2 \,x -\frac{3}{2}+3\, \alpha_2 -\frac{1}{2} \alpha_2^3,$$
$$\alpha_2=\frac{1}{2} \alpha_1+\frac{1}{2}\sqrt{\alpha_1^2+4} \, \mbox{ is a root of } x^2-\alpha_1 \,x-1,$$
 $$\mbox{ and } \alpha_1=\frac{-1}{2}+\frac{1}{2}\sqrt{17} \, \mbox{ is a root of } x^2+x-4\,.$$

 Noticing that $\cos{\frac{2 \pi}{17}}=\frac{\zeta_{17}+\zeta_{17}^{-1}}{2}=\frac{1}{2} \alpha_3,$ Gauss obtains the radical expression 

 \begin{center}
$
\cos \frac {2\pi}{17}=\frac {-1+\sqrt{17}+\sqrt{34-2\sqrt{17}}+2\sqrt{17+3\sqrt{17}-\sqrt{34-2\sqrt{17}}-2\sqrt{34+2\sqrt{17}}}}{16}
,$
\end{center}
  
\noindent  of which he seems to be very proud, for he announces it as the first discovery related to the problem of the division of the circle after two thousand years since the times of Euclid, when the geometric constructions for the regular triangle and pentagon were already known.

\bigskip
At the begining of Chapter VII of the Disquisitiones, Gauss claims the results he is about to show on circular functions are extensible to other transcendental functions, like the ones depending on the elliptic integral of the first kind $$s(r)=\int_0^r \frac{dx}{\sqrt{1-x^4}},$$
which determines the arc length of the lemniscate, the curve 
given by the 
affine model  $$(x^2+y^2)^2=x^2-y^2;$$
alternatively, in polar coordinates by $r^2=\cos(2\theta)$.
This plane curve was introduced by Jakob Bernoulli 
as a multiplicative counterpart of the ellipse. See, for instance,  \cite{RefSch} and \cite{RefSt} 
for more information
 on the lemniscate. Rosen \cite{RefR} gives more details on Gauss's knowledge about the lemniscate and its division with ruler and compass.

\bigskip
Abel took Gauss's suggestion under consideration, obtaining his pleasant result. 
The main ingredient in his proof  of the theorem is the elliptic function $\varphi(s)$, called the Jacobi sine function, that inverts the lemniscate arc-length
function.
The inverse $r = \varphi(s)$ exists and can be 
viewed as an elliptic function
on $\mathbb{C}$ with respect to the lattice 
$
\Lambda = \mathbb{Z} \, (1+i)\,\omega   + \mathbb{Z}\,  (1-i)\,\omega \,,
$
where $$\omega =  2\displaystyle{\int_0^1 \frac{dx}{\sqrt{1-x^4}}} = 
2.622057....$$ is the length of a petal of the lemniscate (see \cite{RefC}). One has $\varphi(\overline{z})=\overline{\varphi(z)}$ where the bar stands for complex conjugation.

\begin{figure}[H]
\begin{center}
\begin{tikzpicture}
    \node[anchor=south west,inner sep=0] (image) at (0,0) 
    {\includegraphics[width=0.7\textwidth]{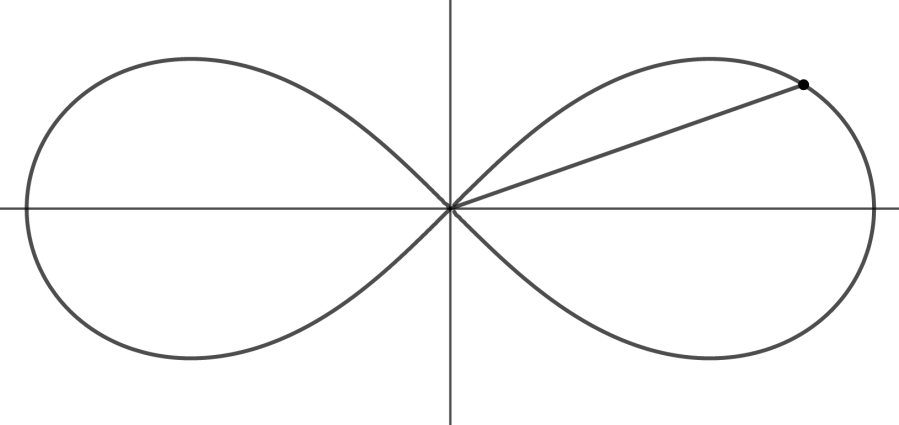}};
    \begin{scope}[x={(image.south east)},y={(image.north west)}]
     \draw[color=black] (0.75,0.65) node {$r$}; 
     \draw[color=black] (0.6,0.8) node {$s(r)$};  
   \end{scope}
\end{tikzpicture}
\end{center}
\caption{The lemniscate arc-length function}
\label{fig:1}       
\end{figure}

The points of division of the lemniscate in $n$ equal parts are determined by the algebraic numbers $\varphi(2\,k\, \omega/n),$ which  play the role of $\zeta_n^k.$ Abel obtains a polynomial with rational coefficients vanishing at $\varphi(2\omega/n)$  after a recursive application of Euler's addition law for the elliptic integrals of the first kind, first proved by Fagnano (1718) in the special case of the lemniscate: 
$$
s(t)= s(r)+s(u) \textrm{ \ if and only if \  }
t = \frac{r\sqrt{1-u^4}+u\sqrt{1-r^4}}{1+r^2u^2}\,.
$$

\bigskip
The purpose of this paper is to give a construction for the regular 17-gon on the lemniscate. 
Instructions to divide the length of the lemniscate in 3 and 5 equal parts with ruler and compass 
 are dealt by the authors 
 in \cite{GoLa}. 
 See also \cite{Osler} for related questions on the lemniscate.
 Before undertaking the division in 17 parts, we devote a section to auxiliary, but very important, constructions on the lemniscate: halving, doubling, adding and
subtracting lemniscate arcs. All these procedures are mainly based on Fagnano's addition formula. 
The last 
section is devoted to the 17-gon. 
The first step in the development of the construction algorithm consists in expressing $\varphi(2\omega/17)$ by radicals. Our first approach to this problem  relied on Gauss's procedure and extended it, computing a tower of real quadratic extensions
$$\mathbb Q(\cos(2\pi/17))\subseteq \mathbb Q(\sin(2\pi/17)) \subseteq \mathbb Q(\sin(2\pi/17), \sqrt{\varepsilon})\subseteq \mathbb Q(\varphi^4(2\omega/17)),$$
where $\varepsilon=1+
\frac{1}{2} \, \sqrt{17} - \frac{1}{2} \, \sqrt{4 \, \sqrt{17} + 17}.$  
Unlike the cyclotomic case, 
the lemniscate extension 
$\mathbb Q(i,\varphi(2\omega/17))/\mathbb Q(i)$ is not cyclic, and one has to be very careful with the choice of the intermediate fields and their primitive elements in order to obtain a radical expression for $\varphi(2\omega/17)$ as simple as possible. However, it turns out that a  better strategy is to apply Fagnano's rule to
$$
\varphi\left(\frac{2\omega}{17}\right) = 
\varphi\left(\frac{\omega}{1+4i}+\frac{\omega}{1-4i}\right),
$$
as Abel did in \cite{Abel2} after extending the addition law to the field of complex numbers in \cite{Abel}. 
Since $1+4i$ is a prime in the ring of gaussian integers, 
the quartic extension $\mathbb Q(i) \subseteq \mathbb Q(\varphi^4(\omega/(1+4i))$ is cyclic (see \cite{RefC}) and $\varphi^4(\omega/(1+4i))$ is easily expressible by radicals, as we show in Section 3. This provides quite a reasonable radical expression for $\varphi(2\omega/17)$ and the corresponding construction with ruler and compass of the points of division of the lemniscate in $17$ equal parts.

\section{Auxiliary constructions}
In this section we shall explain instructions to carry out with ruler and compass the following operations: 
halving, doubling, adding, subtracting, and transferring lemniscate arcs.
Throughout, we assume that the reader is familiarized with the  ruler and compass elementary constructions (for instance, we refer to \cite{RefSu}.)

\begin{figure}[H]
\begin{center}
\begin{tikzpicture}
    \node[anchor=south west,inner sep=0] (image) at (0,0) 
    {\includegraphics[width=0.7\textwidth]{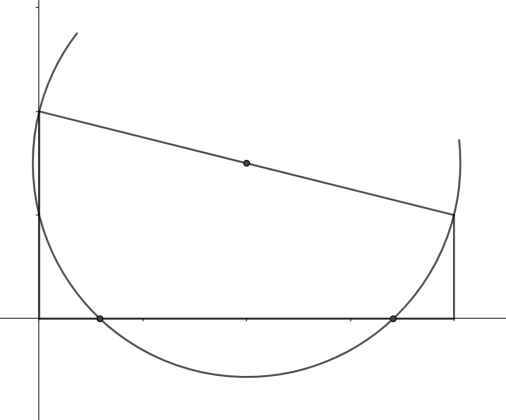}};
    \begin{scope}[x={(image.south east)},y={(image.north west)}]
     \draw[color=black] (0.65,0.75) node {$x^2+2x+1/2$}; 
     \begin{scriptsize}
     \draw[color=black] (0.22,0.12) node {$1-\sqrt{1/2}$};  
     \draw[color=black] (0.77,0.12) node {$1+\sqrt{1/2}$}; 
     \draw[color=black] (0.49,0.2) node {$1$}; 
      \draw[color=black] (0.9,0.2) node {$2$}; 
       \draw[color=black] (0.045,0.49) node {$1/2$}; 
       \draw[color=black] (0.05,0.73) node {$1$}; 
     \end{scriptsize}
   \end{scope}
\end{tikzpicture}
\end{center}
\caption{Negative of the solutions of $x^2+2\,x+1/2=0$}
\label{fig:2}       
\end{figure}
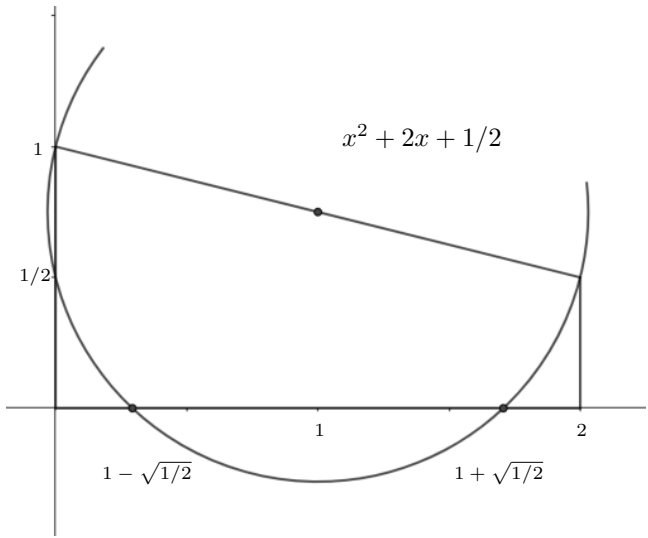

As it is well known, ruler and compass constructions translate into the geometrical construction of roots of linear and quadratic equations. We shall make an intensive use of the
RAT method ({\it right angled trapezium}) to construct the two roots of a quadratic equation at once. Let us recall this method
with an illustrative example: to construct the roots of $x^2+2x+1/2=0$ one forms the right angled trapezium with sides equal to the coefficients lengths $1$, $2$, and $1/2$ and then draw the circle with diameter the transverse side. The two crossing points in the base yield the negative of the two solutions of the quadratic equation (see Fig. \ref{fig:2}).

For the construction in the general case $ax^2+bx+c=0$ 
and some historical background we refer the reader to \cite{CtK}. 
Also, for future reference, we shall denote by $O=(0,0)$ the origin of coordinates, and the points $I=(1,0)$ and $I^*=(-1,0)$. Throughout, the notation $PQ$ will denote the length of the segment with edges the points $P$ and $Q$.

\subsection{Halving arcs.} 
\begin{figure}[H]
\begin{center}
\begin{tikzpicture}
    \node[anchor=south west,inner sep=0] (image) at (0,0) 
    {\includegraphics[width=0.7\textwidth]{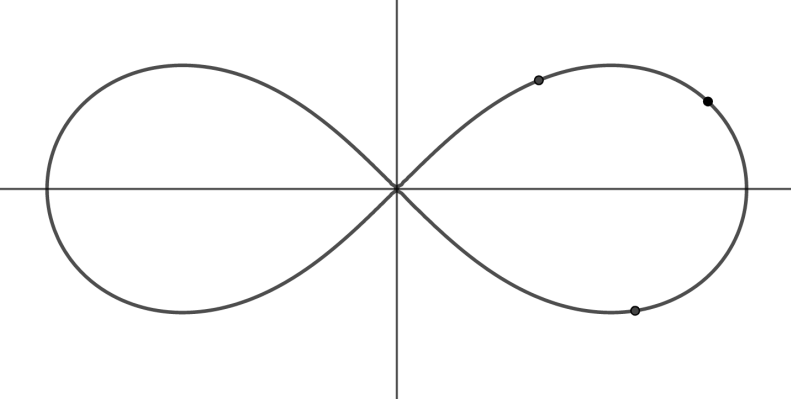}};
    \begin{scope}[x={(image.south east)},y={(image.north west)}]
     \draw[color=black] (0.67,0.85) node {$r_{\theta}$}; 
      \draw[color=black] (0.93,0.755) node {$u_{\varphi}$}; 
     \draw[color=black] (0.8,0.17) node {$r_{\theta}$};  
   \end{scope}
\end{tikzpicture}
\end{center}
\caption{Halving the lemniscate arc determined by $u_\varphi$}
\label{fig:3}     
\end{figure}
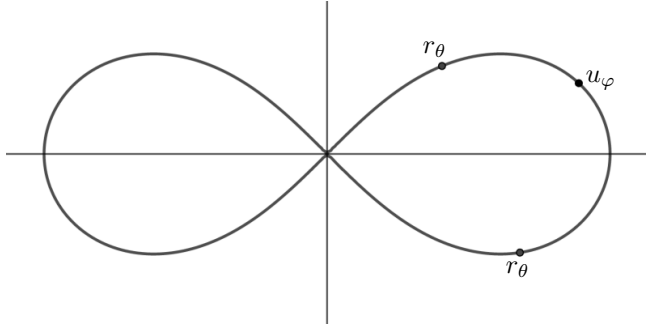

In this section we consider the arc-half practical construction. Given a lemniscate 
point~$u_\varphi$, in polar coordinates, we want to construct with ruler and compass the point~$r_\theta$ such that the arc-length relation 
$2\, s(r) = s(u)$ holds. In fact, in the same process, we shall get two points $r_\theta$ halving each of the two complementary arcs separated 
by~$u_\varphi$ in its petal (see Fig. \ref{fig:3}).

\begin{figure}[H]
\begin{center}
\begin{tikzpicture}
    \node[anchor=south west,inner sep=0] (image) at (0,0) 
    {\includegraphics[width=0.7\textwidth]{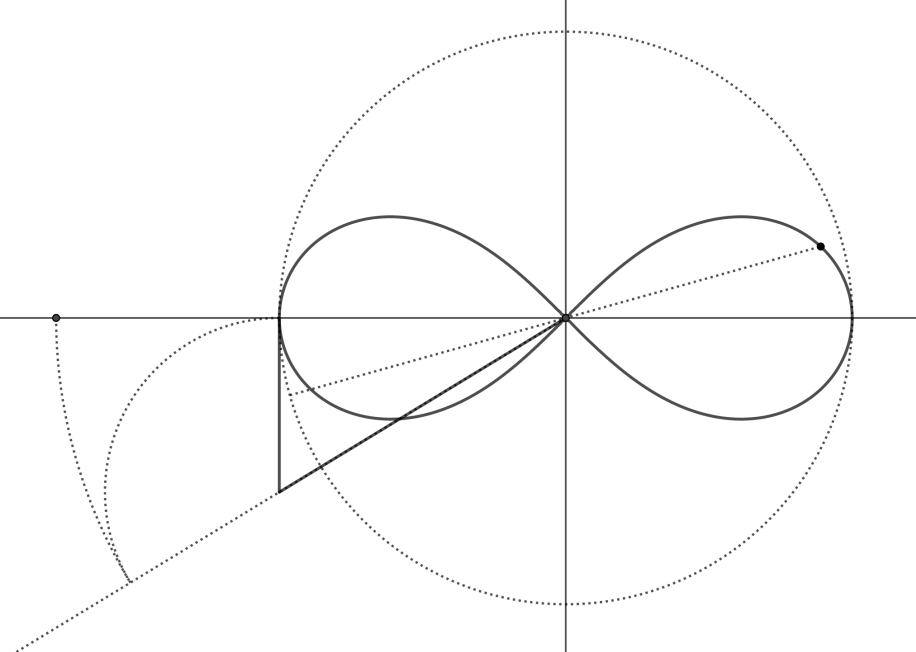}};

    \begin{scope}[x={(image.south east)},y={(image.north west)}]
      \begin{scriptsize}
     \draw[color=black] (0.605,0.46) node {$O$}; 
     \draw[color=black] (0.89,0.66) node {$u_{\varphi}$}; 
     \draw[color=black] (0.45,0.48) node {$\varphi$}; 
     \draw[color=black] (0.45,0.41) node {$\varphi$}; 
     \draw[color=black] (0.065,0.54) node {$A$}; 
     \end{scriptsize}
   \end{scope}
\end{tikzpicture}
\end{center}
\caption{Step 1: Build the point $A$. The segment $OA = \sec(2\varphi)+\tan(2\varphi)$}
\label{fig:4}       
\end{figure}
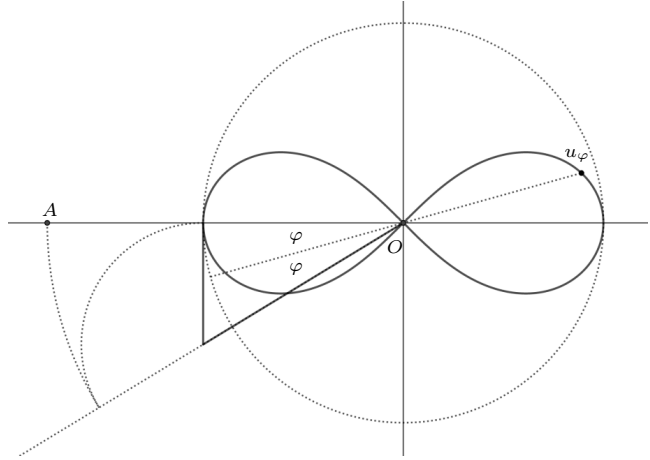

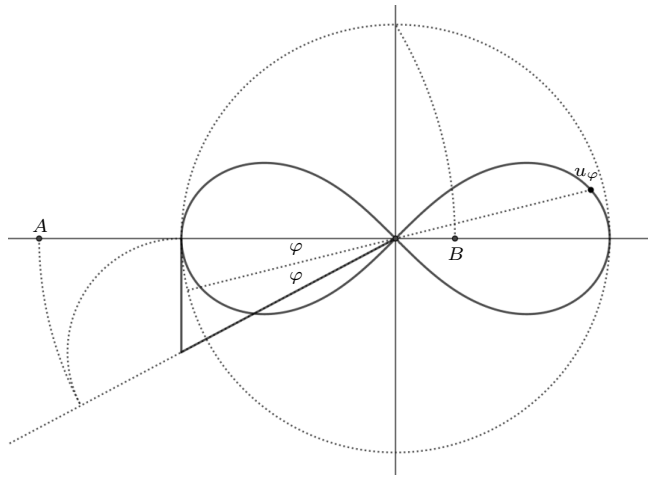
\begin{figure}[H]
\begin{center}
\begin{tikzpicture}
    \node[anchor=south west,inner sep=0] (image) at (0,0) 
    {\includegraphics[width=0.7\textwidth]{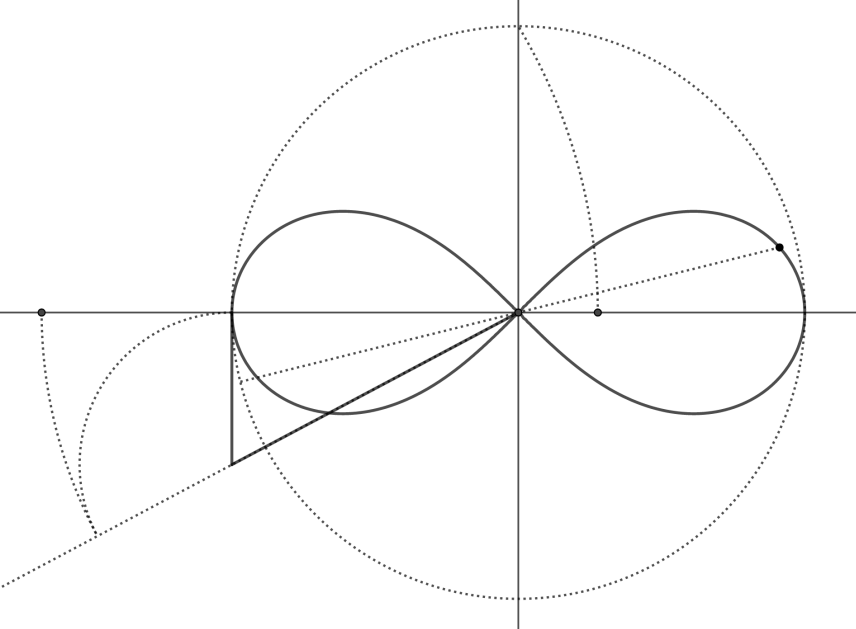}};

    \begin{scope}[x={(image.south east)},y={(image.north west)}]
      \begin{scriptsize}
     \draw[color=black] (0.905,0.64) node {$u_{\varphi}$}; 
     \draw[color=black] (0.45,0.48) node {$\varphi$}; 
     \draw[color=black] (0.45,0.42) node {$\varphi$}; 
     \draw[color=black] (0.05,0.53) node {$A$}; 
       \draw[color=black] (0.7,0.47) node {$B$}; 
     \end{scriptsize}
   \end{scope}
\end{tikzpicture}
\end{center}
\caption{Step 2: Build the point $B$. The segment $OB = \cos(2\theta)$}
\label{fig:5}       
\end{figure}

\begin{figure}[H]
\begin{center}
\begin{tikzpicture}
    \node[anchor=south west,inner sep=0] (image) at (0,0) 
    {\includegraphics[width=0.7\textwidth]{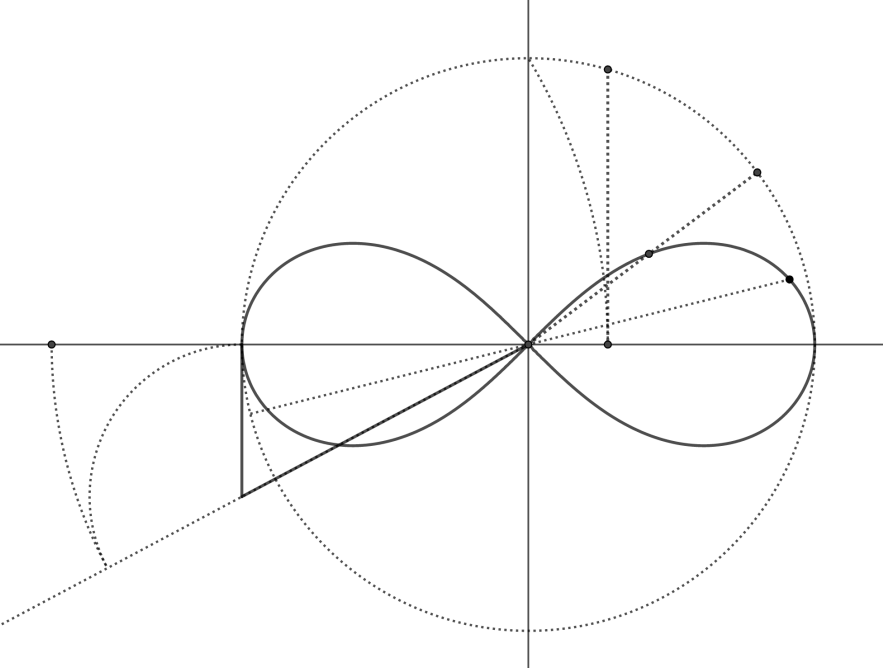}};

    \begin{scope}[x={(image.south east)},y={(image.north west)}]
      \begin{scriptsize}
     \draw[color=black] (0.88,0.62) node {$u_{\varphi}$}; 
     \draw[color=black] (0.45,0.46) node {$\varphi$}; 
     \draw[color=black] (0.45,0.4) node {$\varphi$}; 
     \draw[color=black] (0.06,0.51) node {$A$}; 
       \draw[color=black] (0.685,0.46) node {$B$}; 
         \draw[color=black] (0.7,0.92) node {$C_2$}; 
         \draw[color=black] (0.88,0.76) node {$C_1$}; 
          \draw[color=black] (0.73,0.645) node {$r_{\theta}$}; 
     \end{scriptsize}
   \end{scope}
\end{tikzpicture}
\end{center}
\caption{Step 3: Build the points $C_1$ and $C_2$. The point $r_\theta$ halves the arc corresponding to~$u_\varphi$}
\label{fig:6}      
\end{figure}
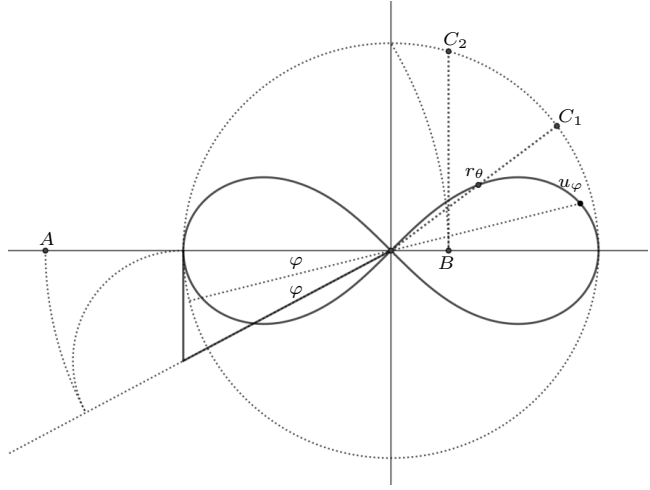

According to Fagnano-Euler formula, the relation between $u$ and $r$ is given by
$$
u = \frac{2r\sqrt{1-r^4}}{1+r^4} \,.
$$
In other words, we need to solve the polynomial equation of degree 8 in $r$:
$$
u^2r^8+4\, r^6+2\, u^2r^4-4\, r^2+u^2 =0\,.
$$
With the variable change $T=r^2=\cos(2\theta)$ and dividing by $u^2$, the above degree 8 polynomial factors as:
$$
\begin{array}{l@{\,} l}
T^4+4/u^2 T^3+2\, T^2-4/u^2 T+1 = & 
\left(T^{2} + 2(\sec(2\varphi)+\tan(2\varphi))T  - 1\right)
\cdot \\[5pt]
& \left(T^{2} + 2(\sec(2\varphi)-\tan(2\varphi))T  - 1\right)\,.
\end{array}
$$
The first quadratic factor yields a bisection point $r_\theta$ of the arc determined by~$u_\varphi$, 
while the second factor provides the point that halves the complementary arc. The ruler and compass instructions to solve the first
factor are as follows (and a similar procedure for the second factor):

\begin{itemize}
    \item[1.] Double the angle $\varphi$ and get the point $A$ such that $OA= 
    \sec(2\varphi)+\tan(2\varphi)$ 
    (see Fig.\,\ref{fig:4}).
    \item[2.] Fold the segment with edges $A$ and $(0,1)$ over the focal axis to get the point $B$ such that $OB = \cos(2\theta)$
    (see Fig.\,\ref{fig:5}).
    \item[3.] Lift a perpendicular to the focal axis through the point $B$ to obtain the point $C_2=1_{2\theta}$ in the unit circle. 
    Bisecting the arc corresponding to $C_2$ one gets $C_1=1_\theta$ and also $r_\theta$ 
    (see Fig.\,\ref{fig:6}).
\end{itemize}

\subsection{Doubling arcs.}

Now we consider the arc-doubling construction.
Given a lemniscate point $r_\theta$, we want to construct the point $u_\varphi$ with ruler and compass in such a way that $s(u) = 2\, s(r)$. To this end,
observe that 
$\sec(2\varphi)+\tan(2\varphi) = \tan \left(\varphi + \frac{\pi}{4}\right)$. 

\begin{figure}[H]
\begin{center}
\begin{tikzpicture}
    \node[anchor=south west,inner sep=0] (image) at (0,0) 
    {\includegraphics[width=0.7\textwidth]{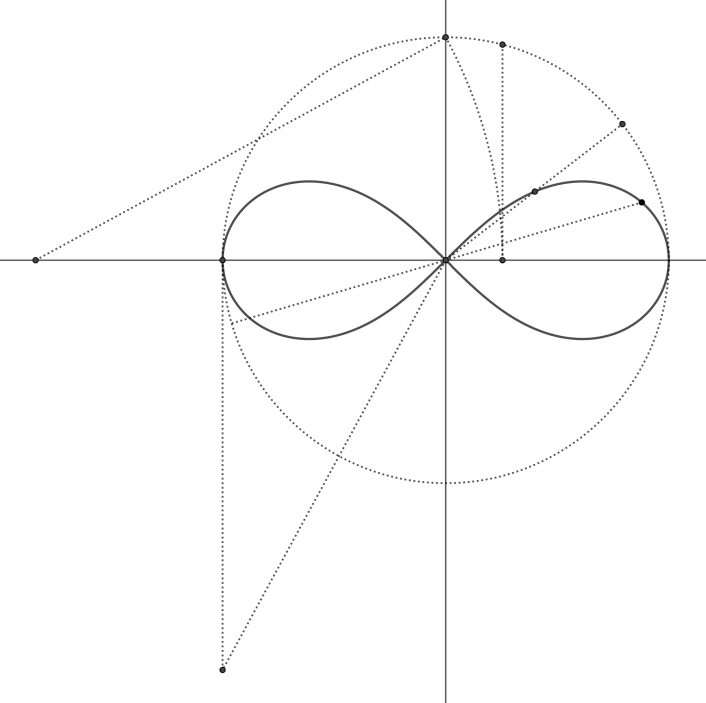}};

    \begin{scope}[x={(image.south east)},y={(image.north west)}]
      \begin{scriptsize}
     \draw[color=black] (0.91,0.73) node {$u_{\varphi}$}; 
     \draw[color=black] (0.75,0.75) node {$r_{\theta}$}; 
     \draw[color=black] (0.45,0.6) node {$\varphi$}; 
     \draw[color=black] (0.5,0.5) node {$\pi/4$};

         \draw[color=black] (0.9,0.84) node {$A$}; 
          \draw[color=black] (0.72,0.96) node {$B$}; 
              \draw[color=black] (0.715,0.61) node {$C$}; 
         \draw[color=black] (0.61,0.965) node {$D$}; 
        \draw[color=black] (0.045,0.61) node {$E$}; 
      \draw[color=black] (0.34,0.045) node {$F$}; 
       \draw[color=black] (0.3,0.61) node {$G$}; 
        
     \end{scriptsize}
   \end{scope}
\end{tikzpicture}
\end{center}
\caption{Doubling a lemniscate arc}
\label{fig:7}       
\end{figure}
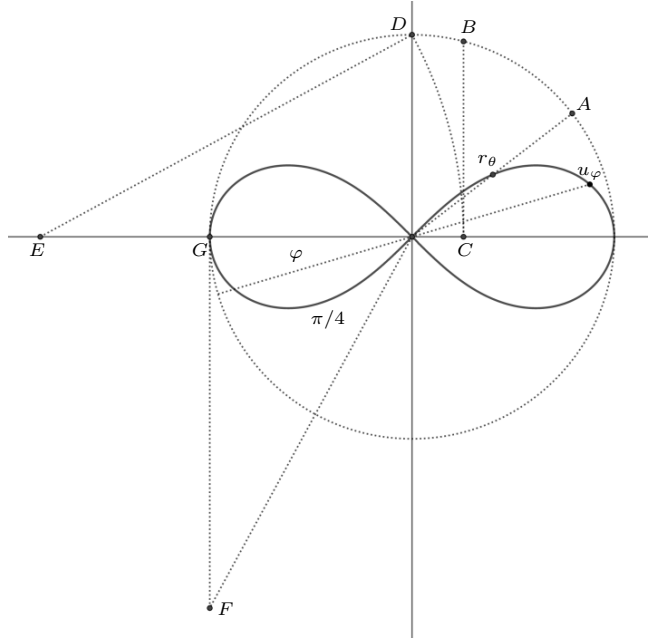

As mentioned in the previous section, the relation between $r_\theta$
and $u_\varphi$ reads:
$$
T^{2} + 2\tan \left(\varphi + \frac{\pi}{4}\right) \, T  - 1 = 0\,,$$
where $T=r^2=\cos(2\theta)$. We can reverse the steps of the above section and get the following instructions (see Fig.\,\ref{fig:7}):

\begin{itemize}
    \item[1.] Double the angle $\theta$ from $A=1_\theta$ to get the point $B=1_{2\theta}$ in the unit 
    circle along with the point $C$ satisfying $OC=\cos(2\theta)$.
    \item[2.] Construct the point $E$ in the focal axis which is the center of a circle passing through $C$ and $D=(0,1)$ or, equivalently, $E$ is the crossing point of the perpendicular bisector of $CD$ with the focal axis. 
    \item[3.] Draw a perpendicular to the focal axis
    through $G=(-1,0)$ to get the point $F$ such that
    $GF=OE=\tan(\varphi +\pi/4)$.
    \item[4.] Subtracting $\pi/4$ to the angle $\angle GOF$, we get the point $u_\varphi$ that doubles the arc.
\end{itemize}

\subsection{Addition and subtraction of arcs.}

Given two points $r_\theta$ and $u_\varphi$ 
in the lemniscate, we want to carry out the addition and subtraction of the corresponding arcs with ruler and compass. That is to say, we need to construct a lemniscate point $t_\beta$ with $s(t)=s(r)+s(u)$, and a lemniscate point $v_\gamma$
such that $s(v)=s(r)-s(u)$.
The addition law implies
$$
t = \frac{r\sqrt{1-u^4}+u\sqrt{1-r^4}}{1+r^2u^2}\,,
\qquad
v = \frac{r\sqrt{1-u^4}-u\sqrt{1-r^4}}{1+r^2u^2}\,.
$$
Equivalently, one checks that $t$ and $v$ are the roots of the quadratic equation
$$
X^2 + B \, X + C = 0
$$
where
$$
B = 
\displaystyle{\frac{-2r\sqrt{1-u^4}}{1+r^2u^2}} \,,
\qquad
C= 
\displaystyle{\frac{r^2-u^2}{1+r^2u^2}}  \,.
$$
In order to apply the RAT method, it is convenient 
to introduce two angles $\alpha$ and $\beta$ such that
$$
\cos(2\theta) = \tan \alpha \,, \quad 
\qquad 
\cos(2\varphi) = \tan \beta \,. 
$$
Then, one has
$$
C = 
\displaystyle{\frac{r^2-u^2}{1+r^2u^2}} = 
\displaystyle{\frac{\cos(2\theta)-\cos(2\varphi)}{1+\cos(2\theta)\,\cos(2\varphi)}} = 
\displaystyle{\frac{\tan \alpha -\tan \beta}
{1+\tan(\alpha)\, \tan(\beta)}} =
\tan(\alpha-\beta) \,.
$$
Similarly, we get the relation
$$
B^2 ={\frac{4\, r^2(1-u^4)}{(1+r^2u^2)^2}} = 
{\frac{4\, \tan\alpha\,  (1-\tan^2\beta)}{(1+\tan\alpha \, \tan \beta)^2}} =
2\, \sin (2\alpha) \, \cos(2\beta) \, \sec^2(\alpha-\beta) \,.
$$
The coefficients $C$ and $B$ can be obtained straightaway from the angles $\alpha$ and $\beta$ as it is displayed in the following diagrams:

\begin{figure}[H]
\begin{center}
  \begin{tikzpicture}
    \node[anchor=south west,inner sep=0] (image) at (0,0) 
    {\includegraphics[width=0.7\textwidth]{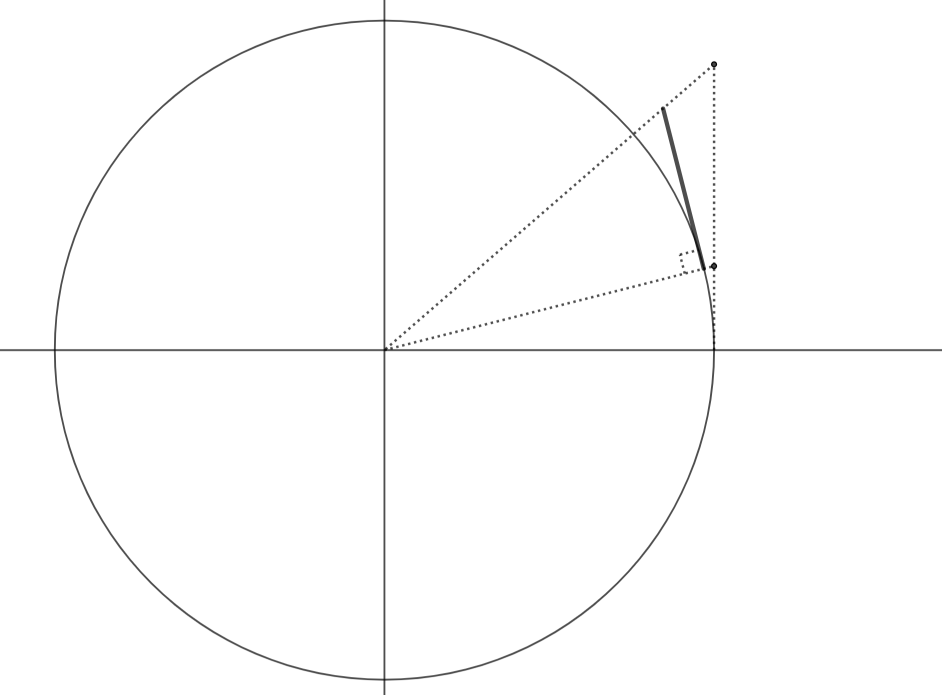}};

    \begin{scope}[x={(image.south east)},y={(image.north west)}]
      \begin{scriptsize}

       \draw[color=black] (0.6,0.53) node {$\beta$}; 
       \draw[color=black] (0.57,0.61) node {$\alpha-\beta$}; 
        \draw[color=black] (0.737,0.76) node {$C$}; 
        \draw[color=black] (0.85,0.91) node {$(1,\cos(2\theta))$}; 
        \draw[color=black] (0.85,0.63) node {$(1,\cos(2\varphi))$}; 
     \end{scriptsize}
   \end{scope}
\end{tikzpicture}
\end{center}
\caption{The coefficient $C=\tan(\alpha-\beta)$}
\label{fig:8}       
\end{figure}
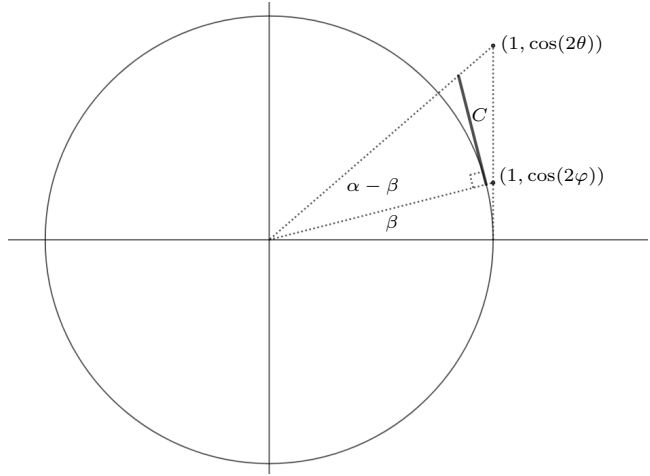

With regard to the coefficient $B$, we use the 
right triangle altitude theorem or geometric mean theorem to construct the length $d=\sqrt{2\, \sin (2\alpha) \, \cos(2\beta)}$.

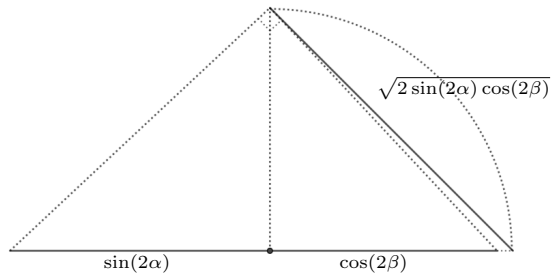
\begin{figure}[H]
\begin{center}
\begin{tikzpicture}
    \node[anchor=south west,inner sep=0] (image) at (0,0) 
    {\includegraphics[width=0.7\textwidth]{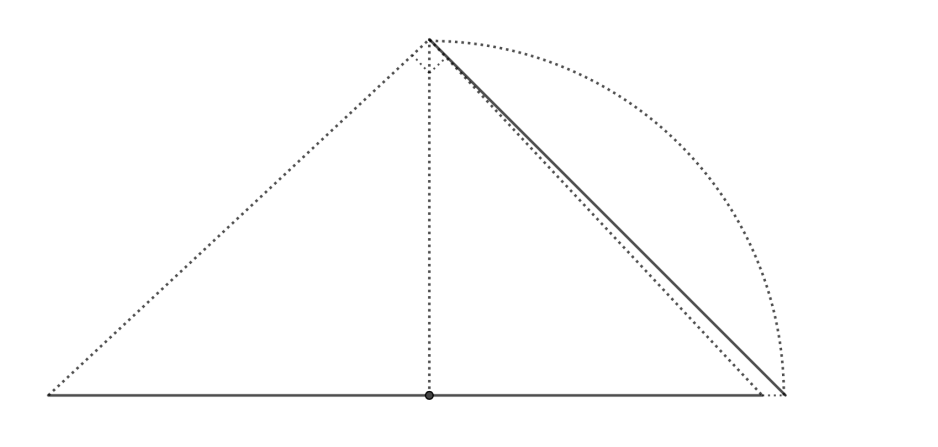}};

    \begin{scope}[x={(image.south east)},y={(image.north west)}]
      \begin{scriptsize}
  
        \draw[color=black] (0.76,0.65) node {$\sqrt{2\sin(2\alpha)\cos(2\beta)}$}; 

      \draw[color=black] (0.25,0.075) node {$\sin(2\alpha)$}; 
    \draw[color=black] (0.62,0.075) node {$\cos(2\beta)$}; 
     \end{scriptsize}
   \end{scope}
\end{tikzpicture}
\end{center}
\caption{The length $d=\sqrt{2\, \sin (2\alpha) \, \cos(2\beta)}$}
\label{fig:9}       
\end{figure}

\begin{figure}[H]
\begin{center}
\begin{tikzpicture}
    \node[anchor=south west,inner sep=0] (image) at (0,0) 
    {\includegraphics[width=0.7\textwidth]{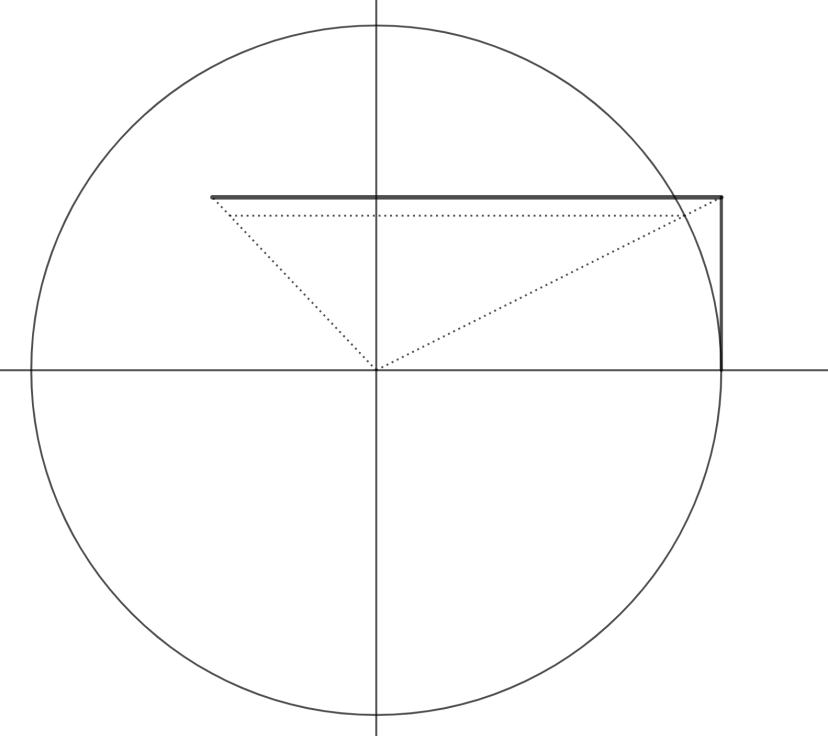}};

    \begin{scope}[x={(image.south east)},y={(image.north west)}]
      \begin{scriptsize}

       \draw[color=black] (0.71,0.62) node {$e$}; 
       \draw[color=black] (0.56,0.68) node {$d$}; 
       \draw[color=black] (0.56,0.755) node {$B$}; 
      
       \draw[color=black] (0.6,0.53) node {$\alpha-\beta$};

        \draw[color=black] (0.89,0.61) node {$C$}; 
     \end{scriptsize}
   \end{scope}
\end{tikzpicture}
\end{center}
\caption{The lengths $e=\sec (\alpha-\beta)$ and 
$B=d\,e$}
\label{fig:10}       
\end{figure}
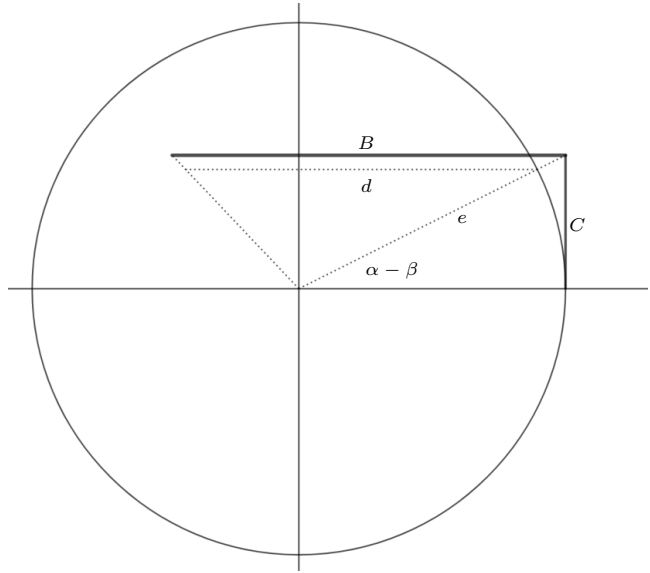

We draw a segment parallel to the $x$-axis 
of length $d$ starting at $1_{\alpha-\beta}$.
Then we use Thales' theorem
to get the segment $B$ drawing a parallel line to the $x$-axis starting 
at $e_{\alpha-\beta}$ of length~$B$, where 
$e=\sec (\alpha-\beta)$ (see Fig.\,\ref{fig:10}).

\begin{figure}[H]
\begin{center}
\begin{tikzpicture}
    \node[anchor=south west,inner sep=0] (image) at (0,0) 
    {\includegraphics[width=0.7\textwidth]{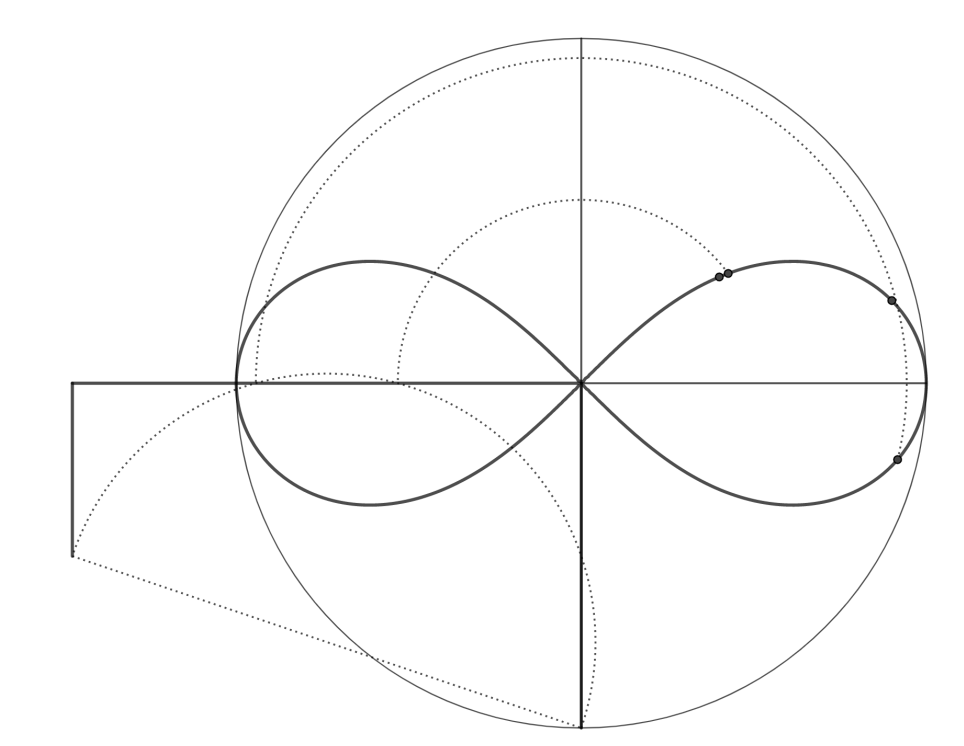}};

    \begin{scope}[x={(image.south east)},y={(image.north west)}]
      \begin{scriptsize}
   
     \draw[color=black] (0.765,0.66) node {$v_{\gamma}$}; 
     
      \draw[color=black] (0.72,0.64) node {$u_{\varphi}$}; 
       
        \draw[color=black] (0.91,0.58) node {$r_{\theta}$}; 

          \draw[color=black] (0.91,0.385) node {$t_{\beta}$};

        \draw[color=black] (0.05,0.38) node {$C$}; 
     
       \draw[color=black] (0.33,0.46) node {$B$}; 
        
     \end{scriptsize}
   \end{scope}
\end{tikzpicture}
\end{center}
\caption{Addition $t_\beta$ and subtraction $v_\gamma$ of the arcs $r_\theta$ and $u_\varphi$}
\label{fig:11}       
\end{figure}
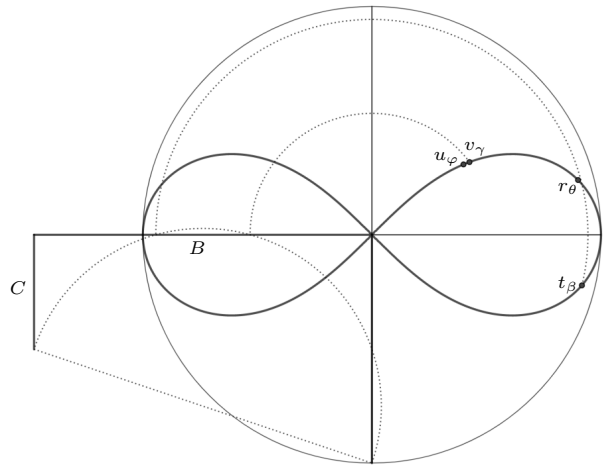

Finally, the right angled trapezium (see Fig.\,\ref{fig:11})
is formed with vertical sides of lengths $-1$ and $-C$ and base of length $-B$. From there, we get 
the addition point $t_\beta$ and the 
subtraction point $v_\gamma$ with ruler and compass.

\subsection{Other constructions.}

The above discussion enables one to practice a wide range of constructions with ruler and compass on the lemniscate. For instance, one can:

\begin{itemize}
    \item[$\bullet$] Transfer an arc to a given point: 
    if the arc is determined by two lemniscate points of radii $r$ and $u$ and the new origin is $t$, then the final edge $w$ is obtained through the operation $s(w) = s(u)-s(r)+s(t)$.
    
    \item[$\bullet$] Bisect an arbitrary arc: let $u_\psi$ be the middle point of the arc given by $r_\theta$ and $t_\varphi$. Then, we have $2(s(u)-s(r))=s(t)-s(r)$. Since we know how to add and halve arcs, we can perform $s(u) = s(r)/2+s(t)/2$ getting the middle point. 
    
\end{itemize}

With regard to $N$-gons, we also have the following specifications:

\begin{itemize}

    \item[$\bullet$]  For $N$ odd, to construct
    the $2N$-gon just draw circles centered in the origin through the vertices of the $N$-gon (see Fig.\,\ref{fig:12}). For $N$ even, drawing the $2N$-gon requires halving only the arcs in the first quadrant.

\begin{figure}[H]
\begin{center}
  \includegraphics[width=0.5\textwidth]{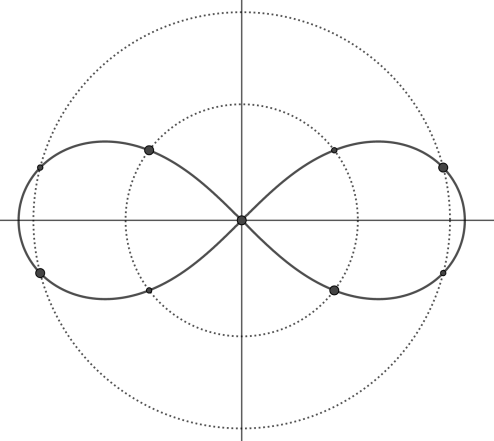}
\end{center}
\caption{The lemniscate regular decagon}
\label{fig:12}    
\end{figure}

    \item[$\bullet$] If the $N$-gon
    and $M$-gon with $\gcd(N,M)=1$ are given, then the B\'ezout identity allows us to find two consecutive vertices (one of each polygon) at distance $2\omega/NM$. Since we know how to transfer arcs to the origin, that gives a procedure to construct the $NM$-gon with ruler and compass.
    
\end{itemize}

\section{The 17-gon} 
  The procedures described in the preceding section, along with Abel's theorem, allow us to construct all regular lemniscate polygons with ruler and compass, as soon as we know how to construct the regular $N$-gon with $N$ being a Fermat prime. The cases $N=3$ and $5$ are dealt in \cite{GoLa}. The case $N=17$ presents some more difficulties, and it is discussed in this section.

\subsection{Radical expressions for
the radii of the lemniscate 17-gon.}

Our first goal  is to discuss 
radical algebraic expressions for
the radii of the lemniscate 17-gon. The points on the lemniscate corresponding to 
$$
r_k = \varphi\left( \frac{2\omega k}{17} \right), \quad \textrm{\  for \ } k=0,\ldots, 16,
$$
form the regular $17$-gon on the lemniscate. 
Recall that $\varphi(s)$ stands for the lemniscate sine elliptic function. As we pointed out in the introduction, we shall follow Abel's idea in \cite{Abel2}, applying the addition law to $r_1:$

$$
\varphi\left(\frac{2\omega}{17}\right) = 
\varphi\left(\frac{\omega}{1+4i}+\frac{\omega}{1-4i}\right)
$$
\begin{equation}\label{radexpr1}=\frac{\varphi\left(\frac{\omega}{1+4i}\right)\sqrt{1-\varphi^4\left(\frac{\omega}{1-4i}\right)}+\varphi\left(\frac{\omega}{1-4i}\right)\sqrt{1-\varphi^4\left(\frac{\omega}{1+4i}\right)}}{1+\varphi^2\left(\frac{\omega}{1+4i}\right)\varphi^2\left(\frac{\omega}{1-4i}\right)}
\end{equation}
where
 $\varphi\left( \frac{\omega}{1-4i}\right)$ is the complex conjugate of $\varphi\left( \frac{\omega}{1+4i}\right).$ Giving an expression by radicals for $\varphi^4\left( \frac{\omega}{1+4i}\right)$ we immediately obtain one for $\varphi\left(\frac{2\omega}{17}\right).$

Abel's procedure to find the minimal polynomial of $\varphi^4\left(  \frac{\omega}{1+4i}\right)$ consists in applying the addition law recursively to $\varphi\left( (1+4i)z\right):$
$$\varphi\left( (1+4i)z\right)=\varphi(z)\frac{P(\varphi^4(z))}{Q(\varphi^4(z))},$$
where
$P(z)=z^{4}+(12-20i)z^{3}-(10-28i)z^{2}-(20+12i)z+1+4i$ and $Q(z)=z^{4}P(1/z).$ Since $\varphi(\omega)=0,$ it turns out that $\varphi^4\left(\frac{\omega}{1+4i}\right)$ is a root of $P(z).$ We know, after a theorem by Eisenstein \cite{Eisenstein}, that this is the minimal polynomial 
of $\varphi^4\left(\frac{\omega}{1+4i}\right)$ over $\mathbb Q(i).$ Its Galois group is cyclic, isomorphic to $\left(\mathbb Z[i]/(1+4i)\mathbb Z[i]\right)^{\times}.$ Notice that $1+4i$ is a prime in $\mathbb Z[i]$ and that $P(z)$ is then the equivalent to the Gauss cyclotomic polynomial. This is an important fact that allows Abel to propose using a Lagrange resolvent (see \cite{Abel, Abel2} or section 12.1.E in \cite{Cox}) to solve this type of equations. That way, 
knowing that the roots of $P(z)$ are 
$$\varphi^4\left(\frac{k\omega}{1+4i}\right), \; \mbox{ with } [k] \in \left(\mathbb Z[i]/(1+4i)\mathbb Z[i]\right)^{\times},$$
we get
$$\varphi^4\left(\frac{\omega}{1+4i}\right)=(-3+5i)-3i\sqrt{1+4i}+(4+i)\sqrt[4]{1+4i}+(-2+i)\sqrt[4]{(1+4i)^3}.$$

\subsection{Construction of the 17-gon by ruler and compass.}

Intending to obtain a construction that takes place mainly
inside the unit circle and trying to avoid situations where two points are indistinguishably close, we have rewritten  the expression above as:

$$\frac{1}{4}\varphi^4\left(\frac{\omega}{1+4i}\right)=
P+Q^{*}+ \frac{3}{2}\sqrt{|Q|}_{\phi/2-\pi/2}  + $$
$$8 \left(|P|_{\theta+3\phi/4}-\frac{1}{2}\sqrt{|Q|}_{3\pi/2-3\phi/4}\right)r^3.
$$
where here $P=-1/2+i/4,\;  Q=1/4+i, \; \theta=\arg(P), \; \phi=\arg(Q) \; \mbox{ and }$
$ Q^*=-1/4+i$ is the reflection of $Q$ across the $y$-axis.  We denote by $|P|$ the modulus of the complex number represented by the point $P$ and by $\arg(P)$ its argument. 
Observe that $\varphi^4\left(\frac{\omega}{1+4i}\right)$ is constructible with ruler and compass using only the points $P$ and $Q.$

\bigskip
The classical procedure for constructing the square root of a positive real number is a particular case of the RAT method used in Section 2 (see, for instance, section 1.3 \cite{Carrega}). Also, multiplication of positive real numbers is based on Thales's theorem, whereas addition and subtraction of two complex numbers  follow the parallelogram law for vectors. 
With that in mind, let's describe a method for constructing 
$\varphi^4\left(\frac{\omega}{1+4i}\right)$ with ruler and compass:

\begin{enumerate}
\item Mark $O=(0,0), I=(1,0), J=(0,1), P, Q, Q^{*}$, and the corresponding angles $\theta$ and $\phi.$
\item Construct a segment of length $s=\frac{1}{2}\sqrt{|Q|}$ and then its square root $r=\sqrt{s}.$ 
\item Mark the points $S=(0,s)$ and $R=(0,r).$
\item Draw circles $C_1$ and $C_2$ centered at $O$ and passing through $P$ and $S$ respectively.
\item Construct the bisector of the angle $\angle{IOQ},$ and mark its point of intersection $T$ with $C_1.$ Then, draw the bisector of $\angle{TOQ}.$ 
We have determined the angles $\phi/2$ and $3\phi/4.$ 
\item Mark $A=|P|_{\theta+3\phi/4}$ on $C_1,$ and $B=|S|_{3\pi/2-3\phi/4}, C=|S|_{\phi/2-\pi/2}$ on $C_2.$
\item Translate the segment $\overline{AB}$ to the parallel line passing through $O,$ taking $B$ to $O$ and $A$ to $D=A-B.$ Mark $E=8D.$
\item The parallel line to the segment of edges $J, E$ passing through $R$ intersects the ray $\overline OE$ at $Er,$ whereas the parallel line to the segment of edges $J, Er$ passing through $S$ intersects the same ray at $Er^3.$
%
\item  Mark $U=\varphi^4\left(\frac{\omega}{1+4i}\right)=4\left((Er^3+C)+(P+C)+(Q^*+C)\right).$

\end{enumerate}

\begin{remark} 
The addition $(Er^3+C)+(P+C)+(Q^*+C)$ provides such a short segment that the angle formed by $(Er^3+C)+(P+C),$ $O$ and $Q^*+C$ is close to $\pi$ and thus the determination by eyesight of the fourth point of the parallellogram is difficult. It is much easier to determine the resulting point of the subtraction
$(Er^3+C)+(P+C)-(-(Q^*+C)),$ as we have done with points $A$ and $B$ before.
\end{remark}

\begin{figure}[H]
\begin{center}
 \begin{tikzpicture}
    \node[anchor=south west,inner sep=0] (image) at (0,0) 
    {\includegraphics[width=\textwidth]{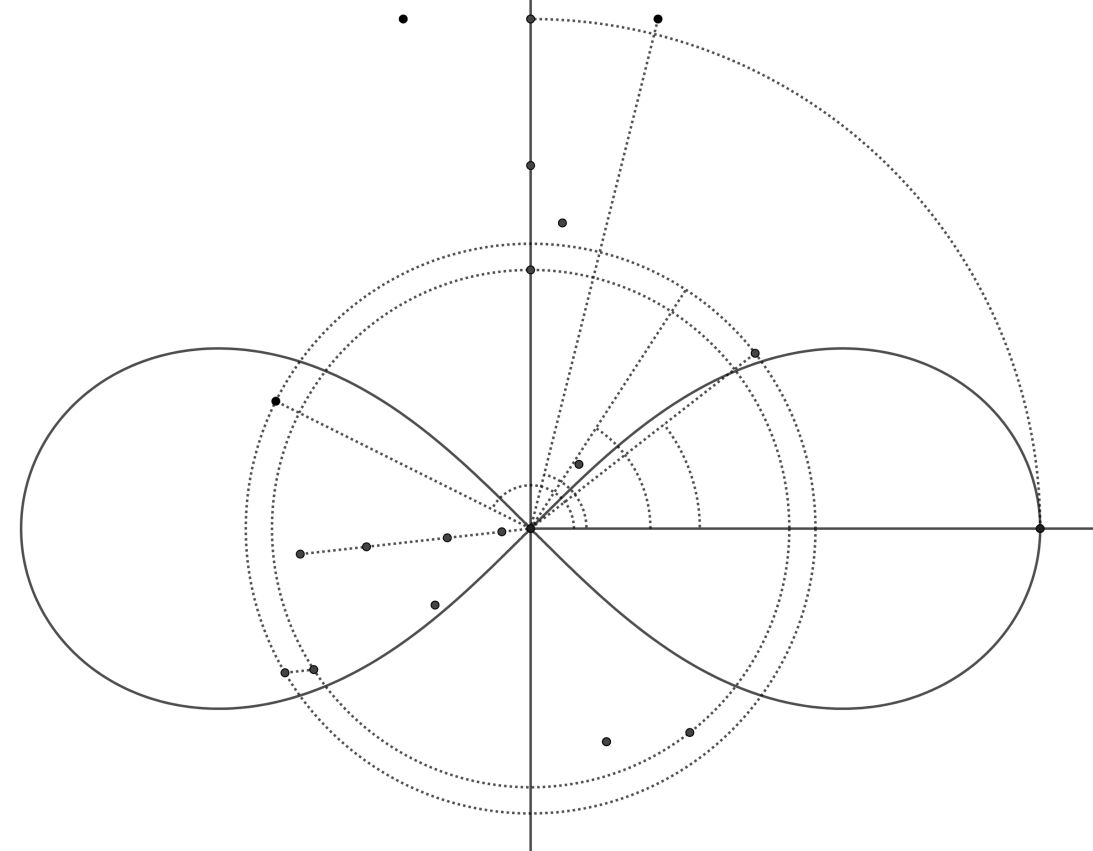}};
    \begin{scope}[x={(image.south east)},y={(image.north west)}]
\begin{scriptsize}

\draw[color=black] (0.62,0.155) node {$C$}; 

\draw[color=black] (0.53,0.15) node {$Er^3\!+\!C$}; 

 \draw[color=black] (0.62,0.975) node {$Q$};

\draw[color=black] (0.39,0.975) node {$Q^{*}$}; 

 \draw[color=black] (0.96,0.393) node {$I$}; 

 \draw[color=black] (0.47,0.975) node {$J$}; 

 \draw[color=black] (0.47,0.81) node {$R$}; 

  \draw[color=black] (0.47,0.665) node {$S$}; 

\draw[color=black] (0.525,0.76) node {$Q^{*}\!\!+\!C$};

\draw[color=black] (0.7,0.6) node {$T$};

\draw[color=black] (0.51,0.367) node {$O$}; 

\draw[color=black] (0.41,0.39) node {$Er^3$}; 
\draw[color=black] (0.36,0.29) node {$P\!+\!C$};
\draw[color=black] (0.33,0.379) node {$Er$}; 
\draw[color=black] (0.275,0.37) node {$E$}; 
\draw[color=black] (0.455,0.358) node {$D$}; 

\draw[color=black] (0.465,0.44) node {$\theta$};

\draw[color=black] (0.52,0.47) node {$U$}; 

\draw[color=black] (0.55,0.4) node {$\phi$}; 

  \draw[color=black] (0.61,0.41) node {$\frac{3\phi}{4}$};     

  \draw[color=black] (0.655,0.42) node {$\frac{\phi}{2}$};     
        
      \draw[color=black] (0.245,0.21) node {$A$};  
       \draw[color=black] (0.303,0.22) node {$B$};  
        \draw[color=black] (0.24,0.54) node {$P$};

     \end{scriptsize}
   \end{scope}
\end{tikzpicture}
\end{center}
\caption{$U=\varphi^4\left(\frac{\omega}{1+4i}\right)$}
\label{fig:13}    
\end{figure}
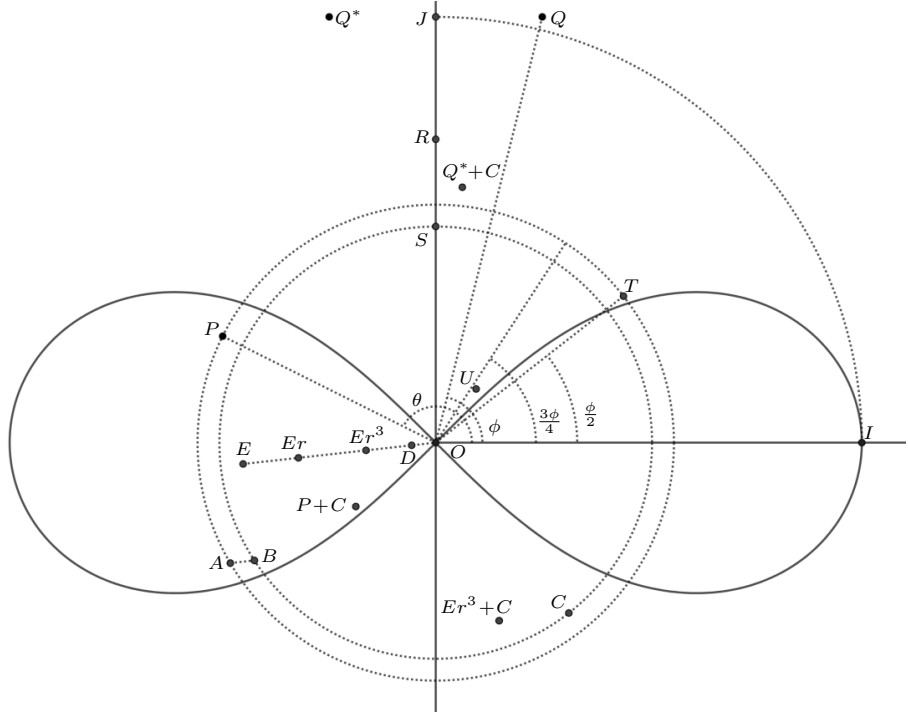

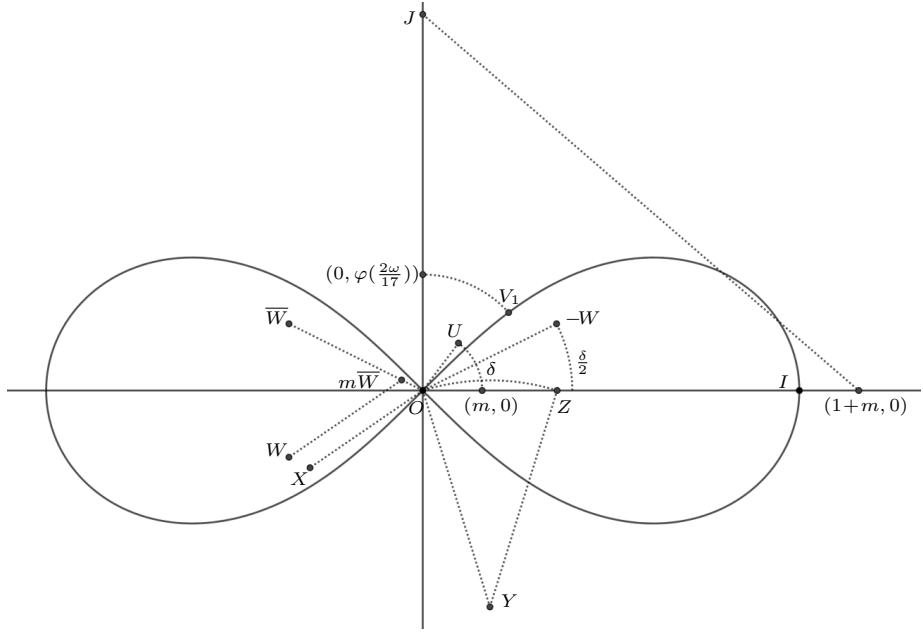
\begin{figure}[H]
\begin{center}
 \begin{tikzpicture}
    \node[anchor=south west,inner sep=0] (image) at (0,0) 
    {\includegraphics[width=\textwidth]{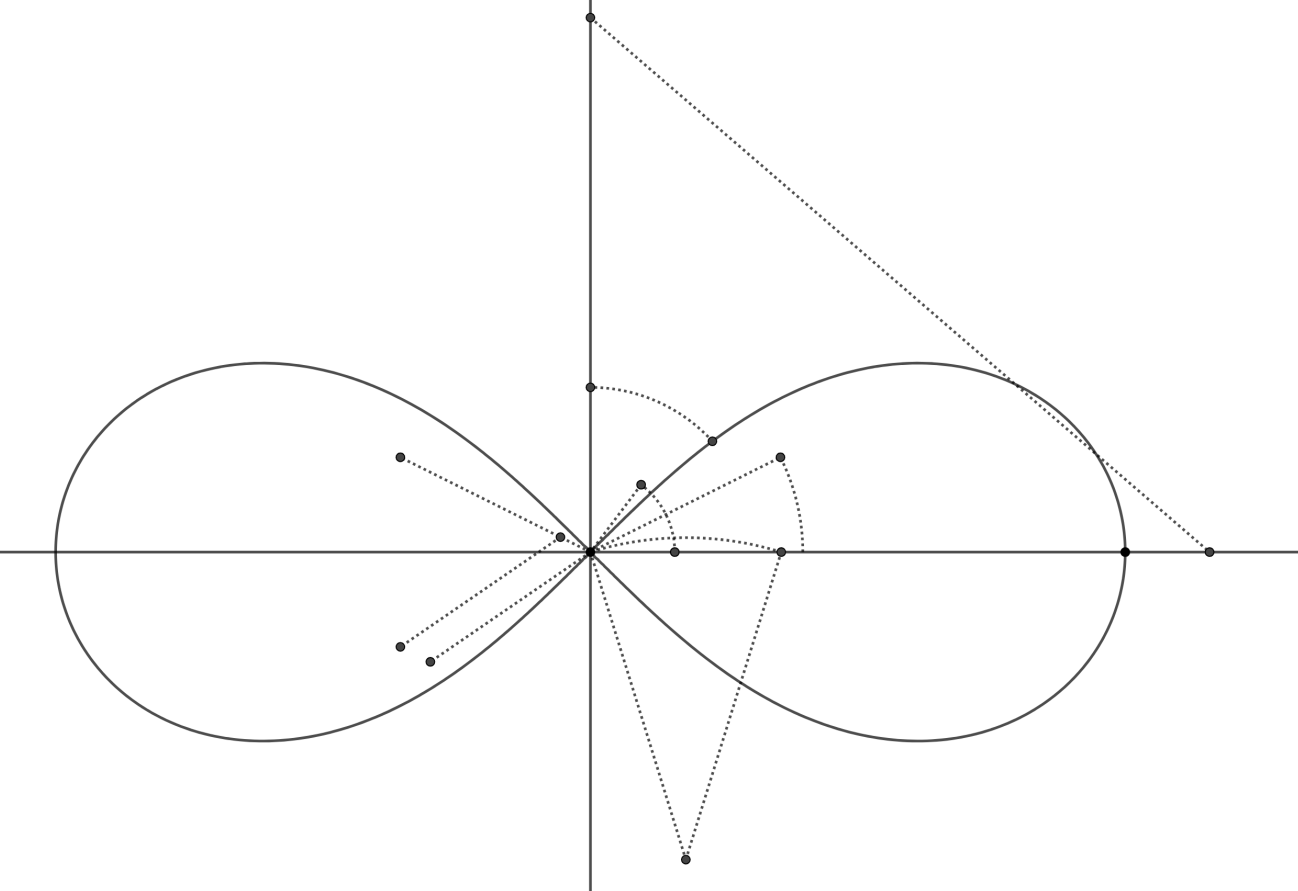}};
    \begin{scope}[x={(image.south east)},y={(image.north west)}]
\begin{scriptsize}

\draw[color=black] (0.55,0.045) node {$Y$}; 

 \draw[color=black] (0.94,0.355) node {$(1\!+\!m,0)$}; 
 \draw[color=black] (0.85,0.395) node {$I$}; 

 \draw[color=black] (0.44,0.975) node {$J$}; 
 
  \draw[color=black] (0.4,0.565) node {$(0,\varphi(\frac{2\omega}{17}))$}; 

\draw[color=black] (0.295,0.285) node {$W$};
\draw[color=black] (0.32,0.24) node {$X$};
\draw[color=black] (0.295,0.5) node {$\overline W$}; 

  \draw[color=black] (0.63,0.495) node {$-\!W$};    

  \draw[color=black] (0.55,0.53) node {$V_1$};   

\draw[color=black] (0.448,0.352) node {$O$}; 

\draw[color=black] (0.385,0.398) node {$m\overline W$};

\draw[color=black] (0.49,0.475) node {$U$}; 

\draw[color=black] (0.53,0.415) node {$\delta$}; 

\draw[color=black] (0.53,0.355) node {$(m,0)$}; 
\draw[color=black] (0.61,0.355) node {$Z$};

  \draw[color=black] (0.63,0.425) node {$\frac{\delta}{2}$};

     \end{scriptsize}
   \end{scope}
\end{tikzpicture}
\end{center}
\caption{First vertex $V_1$ of the lemniscate $17$-gon}
\label{fig:14}    
\end{figure}

We denote by $m$ and $\delta$ the modulus and argument respectively of the point $U=\varphi^4\left(\frac{\omega}{1+4i}\right).$
The complex number $W=\varphi^2\left(\frac{\omega}{1+4i}\right)$ is one of its square roots,
$W=-\sqrt{m}_{\delta/2},$
which is easily drawable by bisecting $\delta$ and constructing the square root of the length of the segment $OU$ (see Figure \ref{fig:14}).
Rewriting equation (\ref{radexpr1}) one has

$$\varphi\left(\frac{2\omega}{17}\right)=\frac{2Re\left( \sqrt{W-m\overline W}\right)}{1+m},$$
so that the point $W$ and the radius $m$ are enough to obtain the first vertex of the $17$-gon on the lemniscate. A simple method is described below:
\begin{enumerate}
\item Once $W$ is marked, reflect it over the $x$-axis to get its complex conjugate, $\overline W.$ 

\item  The parallel line to the segment of edges $\overline W, I$ passing through $(m,0)$ intersects the ray $O\overline W$ at $m\overline W.$

\item Translate the segment of edges $W, m\overline W$  to the parallel line passing through $O,$ taking $m \overline W$ to $O$ and $W$ to $X=W-m\overline W.$

\item Construct a segment of length $\sqrt{|X|}$ and bisect the angle $\angle{XOI}$ to get the point $Y=\sqrt{W-m\overline W}.$

\item The circle centered at $Y$ passing through $O$ cuts the semiaxis $x>0$ at the point $Z=\left(2Re\left( \sqrt{W-m\overline W} \right),0\right).$

\item The parallel line to the segment of edges $J, (1+m,0)$ passing through $Z$ intersects the semiaxis $y>0$ at $(0,\varphi(2\omega/17).$ 

\item Fold the segment of edges $O, (0,\varphi(2\omega/17)$ over the lemniscate to get $V_1,$ the first point of division of the lemniscate in $17$ equal parts.  
\end{enumerate}

Vertex $V_2$ is obtained doubling the arc $OV_1,$ whereas $V_3$ is the addition of the arcs based on $O$ and determined by $V_1$ and $V_2.$ Vertices $V_4$ and $V_6$ are obtained doubling the arcs determined respectively by $V_2$ and $V_3.$ Bisecting the arc $OV_1$ and reflecting across the $x$-axis we find $V_8.$ The same operation on the arc $OV_3$ returns $V_7.$ To obtain $V_5$ we can reflect $V_6$ across the $x$-axis and double the corresponding arc based on $O.$ The remainder points of division of the lemniscate in 17 parts are the reflections of $V_1, \dots, V_8$ about $O.$

\begin{figure}[H]
\begin{center}
\begin{tikzpicture}
    \node[anchor=south west,inner sep=0] (image) at (0,0) 
    {\includegraphics[width=\textwidth]{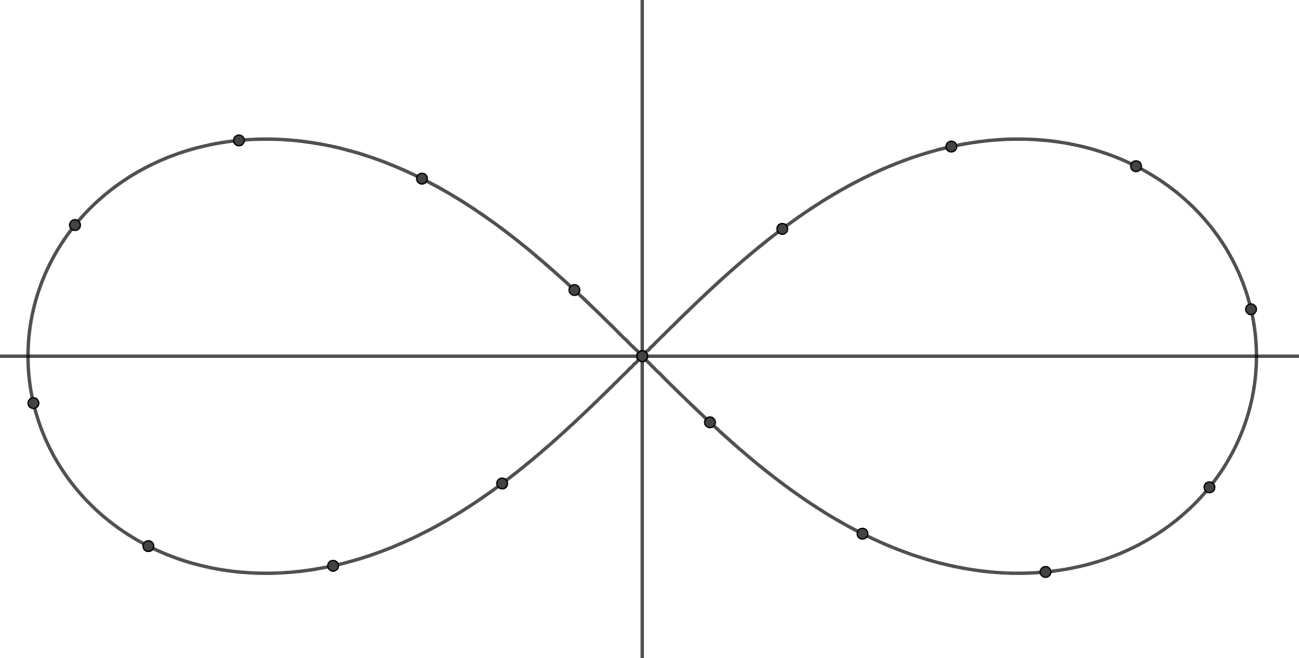}};
    \begin{scope}[x={(image.south east)},y={(image.north west)}]
\begin{scriptsize}

  \draw[color=black] (0.484,0.415) node {$O$};   

\draw[color=black] (0.53,0.34) node {$V_8$}; 

  \draw[color=black] (0.64,0.18) node {$V_7$};     

\draw[color=black] (0.58,0.66) node {$V_1$};   
   
   \draw[color=black] (0.73,0.81) node {$V_2$};       
       
 \draw[color=black] (0.79,0.16) node {$V_6$};   

 \draw[color=black] (0.89,0.77) node {$V_3$}; 
 
\draw[color=black] (0.94,0.53) node {$V_4$};

\draw[color=black] (0.91,0.26) node {$V_5$};
      
     \end{scriptsize}
   \end{scope}
\end{tikzpicture}
 \end{center}
  \caption{All vertices of the 17-gon}
  \label{fig:20}
\end{figure}
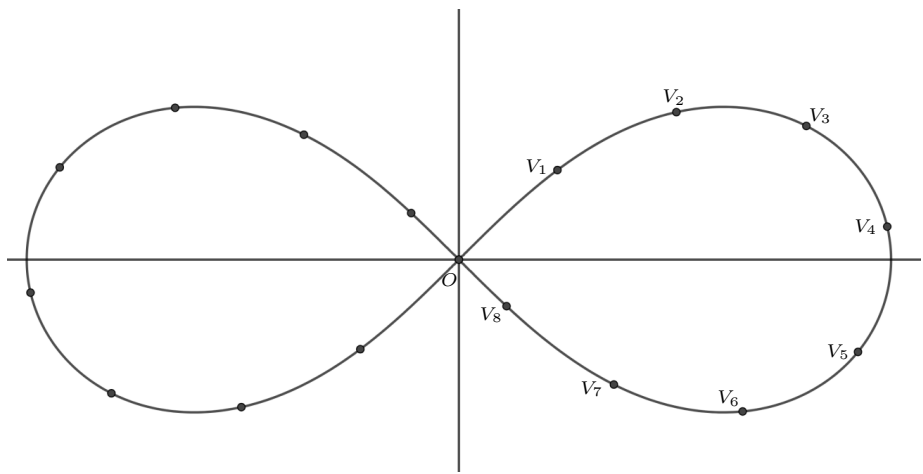

\end{document}